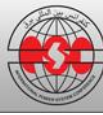

# An Efficient Scenario-Based Stochastic Model for Dynamic Operational Scheduling of Community Microgrids with High Penetration Renewables

Farhad Samadi Gazijahani[1*], Javad Salehi[2]

[1,2]Azarbaijan Shahid Madani University
Department of Electrical Engineering
Tabriz, Iran

[1]f.samadi@azaruniv.ac.ir, [2]j.salehi@azaruniv.ac.ir

*Abstract*—The supply of electrical energy is being increasingly sourced from renewable generation resources. The variability and uncertainty of renewable generation, compared to a dispatchable plant, is a significant dissimilarity of concern to the traditionally reliable and robust distribution systems. In order to reach the optimal operation of community Microgrids (MGs) including various Distributed Energy Resource (DER), the stochastic nature of renewable generation should be considered in the decision-making process. To this end, this paper proposes a stochastic scenario based model for optimal dynamic energy management of MGs with the goal of cost and emission minimization as well as reliability maximization. In the proposed model, the uncertainties of load consumption and also, the available output power of wind and photovoltaic units are modeled by a scenario-based stochastic programming. Through this method, the inherent stochastic nature of the proposed problem is released and the problem is decomposed into a deterministic problem. Finally, an improved metaheuristic algorithm based on Cuckoo Optimization Algorithm (COA) is implemented to yield the best global optimal solution. The proposed framework is applied in the typical grid-connected MGs in order to verify its efficiency and feasibility.

## Nomenclature

| | |
|---|---|
| $s$ | Set of scenarios |
| $b$ | Branch |
| $P_{WT(v)}$ | Generated power at wind turbine |
| $P_{r,WT}$ | Rated power of WT |
| $P_D$ | Load demand |
| $v$ | Wind speed |
| $P_{PV}$ | Generated power by PV unit |
| $G_{ING}$ | Incident brilliance |
| $G_{STC}$ | Brilliance at standard situation |
| $P_{STC}$ | permissive output of PV unit |
| $Z$ | Objective function |
| $Tc, Tr$ | Cell reference temperatures respectively |
| $k$ | Maximum power correction temperature |
| $\alpha, \beta$ | Shape and scale parameter of Beta |
| $c$ | Scale of the Rayleigh model |
| $\Theta, \varrho$ and $\gamma$ | CHPs coefficients |
| $\lambda$ | Failure rate |
| $SOC$ | State of charge of ESS units |
| $\sigma$ | Standard deviation of load demand |
| $P_{ch}/P_{dis}$ | Power charge and discharge of ESS |
| $\zeta, \delta$ | Generator emission characteristics |
| $ENS$ | Energy not supplied |
| $\eta$ | Efficiency of CHP units |
| $\mu$ | Mean value of load demand |
| $C_{int}$ | Price of energy not supply |
| $C_{Ploss}$ | Price of power losses |
| $N_{res}$ | Disjunct nodes within defect |
| $N_{rep}$ | Conjunct nodes within shaving |
| $P_{res}$ | Not remade demand pending error |
| $P_{rep}$ | Remade demand pending error |
| $T_{res}$ | Period of error position |
| $T_{rep}$ | Period of renovation time |
| $R_b$ | Resistance of the line b |
| $L_b$ | Length of the line b |
| $H$ | Weight multipliers |

## I. Introduction

Microgrids (MGs) are power systems comprising distributed energy resources (DERs) and electricity end-users, possibly with controllable elastic loads, all deployed across a limited geographic area [1]. Depending on their origin, DERs can come either from distributed generation (DG) or from Energy Etorage Systems (ESS). DG refers to small-scale power generators such as diesel generators, fuel cells, and renewable energy sources (RES), as in wind or photovoltaic (PV) generation. Since, MGs can participate in power markets and also provide some ancillary services, proper scheduling of the MG is essential from the main grid point of view [2]-[6]. Therefore, a suitable strategy should be pursued for the MG operation.

In recent years, many researchers have gravitated to optimization-based fields [8]-[10]. Locating global optima in many real-world optimization problems is often painstakingly tedious. Thus, finding a better optimization algorithm is a critical task. The self-adaptive evolutionary programming method is implemented for solving the nonlinear optimal power flow problem in [8]. The bacterial foraging algorithm is used for the optimal control of a DG system in [9] and [10].

Stochastic programming is also used to cope with the variability of RES. Single-period chance-constrained economic dispatch problems for RES have been studied in [11], yielding probabilistic guarantees that the load will be served. Considering the uncertainties of demand profiles and PV generation, a stochastic program is formulated to minimize the overall cost of electricity and natural gas for a building in [12]. Without DSM, robust scheduling problems with penalty-based costs for uncertain supply and demand have been investigated in [13]. Recent works explore energy scheduling with DSM and RES using only centralized algorithms [14], [15]. An energy source control for optimal operation of MG is formulated and solved using model predictive control in [16].

The probabilistic methods can be classified into three categories: The Monte Carlo Simulation (MCS), the analytical techniques and the approximate methods. The MCS methods are the most straightforward and accurate one but have the shortcoming of remarkable computational efforts. Analytical techniques need fewer number of simulations but still require complicated mathematical computation. Approximate methods provide balance among computational efficiency and accuracy. The 2m point estimate method (PEM), as approximate method, is an efficient and reliable method to model the uncertainty in power systems.

We focus on selecting a strong and robust evolutionary algorithm in order to solve the proposed complex operational scheduling of MGs that is located on the top of the dispatching priority. One of the new evolutionary algorithms that has great potential is a cuckoo optimization algorithm (COA). Also, this paper implements the stochastic scenario based approach to model the uncertainty in hourly load demand and available output power of solar and wind DG units.

The rest of this paper is structured as follows: Section II shows the stochastic scenario based model. In Section III, the modelling of load and renewable generation are shown. The proposed problem is formulated in Section IV and Section V explains the utilized algorithm; results and conclusions are drawn in Section VI and VII, respectively.

## II. STOCHASTIC ANALYSIS APPROACH

Random sampling is a corner stone of the stochastic scenario based method, hence using a most efficient method for scenario generation is crucial. Load uncertainty roots from the characteristics of the load that can be affected by many different variables such as weather, temperature, humidity, special programs followed by the governments, etc. Here we assume that a prediction tool is available for hour-by-hour sequence of load forecasts. Since the load forecasts are generally inaccurate, the forecast error is modelled as zero-mean normally distributed random variable. Note that the normality assumption of the demand forecast error is standard in the literature [17]. In this regard, firstly the probability density function (PDF) is performed by the use of data samples of previous several years. Then depending on the desired accuracy, the PDF is divided into different levels. The area related to each level represents the probability of such forecast error level. Wind energy as one the most important kinds of the renewable energy sources has played a significant role in the growth of the renewable resources. Owing to its characteristics of the environmental-compatibility and also the economical fuel costs, it has experienced fast development in the recent years. However, as the random variation of the wind speed, the output power generated by the WT is neither continuous nor stable and therefore cannot be supposed as a reliable power source [17]. The basic idea to model these variations in the evaluations is the probabilistic methods.

## III. SYSTEM MODEL

### A. Modeling of DG Units

The model describing each DG is described as follows.

*a) WT*: The power generated by WT as a function of wind speed can be calculated by (1) and modelled by Weibull PDF:

$$P_{WT} = \begin{cases} 0 & 0 < V < P_{STG} \\ (A.V^2 + B.V + C) * P_{rate} & V_{ci} \leq V \leq V_r \\ P_{rate} & V_r \leq V \leq V_{co} \\ 0 & V_{co} \leq V \leq \infty \end{cases} \quad (1)$$

The relationship between the generated power and wind speed of the WT can be demonstrated as Fig.1.

*b) PV*: The power generated by PV depends on the irradiance and the ambient temperature and modelled by Beta PDF:

$$f(s) = \begin{cases} \frac{\Gamma(\alpha+\beta)}{\Gamma(\alpha)+\Gamma(\beta)} \times s^{(\alpha-1)} . (1-s)^{\beta-1} & , 0 \leq s \leq 1, \alpha \geq 0, \beta \geq 0 \\ 0 & , otherwise \end{cases} \quad (2)$$

$$\beta = (1-\mu) . \left(\frac{\mu.(1+\mu)}{\sigma^2} - 1\right) \quad (3)$$

$$\alpha = \frac{\mu.\beta}{1-\mu} \quad (4)$$

$$P_{pv} = P_{STG} * \frac{G_{ING}}{G_{STG}} * (1 + k(T_C - T_{ref})) \quad (5)$$

*c) CHP*: The output power generated by CHP is controlled by an installed governor. The rate of fuel consumed by CHP can be expressed as below.

$$L_{CHP}(P_{CHP}) = \Theta.P_{CHP}^2 + \varrho.P_{CHP} + \gamma \quad (6)$$

*d) ESS*: Since there are some restrictions on charge and discharge rate of storage devices during each time interval, the following equation can be considered:

$$SOC_t = SOC_{t-1} + \eta_{ch} P_{ch} \Delta t - \frac{1}{\eta_{dis}} P_{dis} \Delta t \quad (7)$$

### B. Load Demand Model

The electric load is an uncertain parameter in deregulated environment such as MGs but this parameter is specifically tied with various factors. An increase/decrease in electric load will tend to increase/decrease in electricity price and vice versa. Without loss of generality, the correlation between wind speed and load-price pattern are assumed to be independent which can be modeled as a normal distribution function by mean value ($\mu$) and standard deviation ($\sigma$) as:

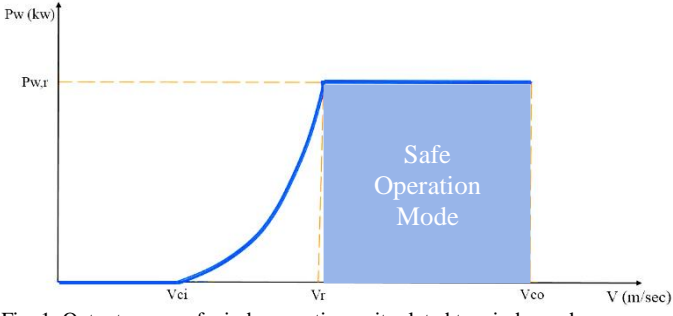

Fig. 1. Output power of wind generation unit related to wind speed.

$$f(p_l) = \frac{1}{\sqrt{2\pi} \times \sigma} exp(-\frac{(p_l - \mu)^2}{2 \times \sigma^2}) \quad (8)$$

## IV. PROPOSED MODEL

### A. Cost evaluation

The objective function of the DED problem is to minimize the total production cost over the operating horizon by minimizing $F1$. The total energy and operating cost of the MG includes the fuel costs of units, operation costs, emission cost as well as power losses costs. The utilized objective function, simultaneously, minimizes fuel cost and operating costs. The optimal values, obtained by minimizing $F1$, provide the amount of power generated by each DG, as well as power sold to or purchased from the main grid in day ahead scheduling.

$$F1 = \sum_{t=1}^{T}\sum_{s=1}^{S}\sum_{n=1}^{N}[f_{n,t}^s + OM_{n,t}^s + E_{n,t}^s + PLC_{loss}] \quad (9)$$

*a) Fuel Consumption Function*: The operating cost in fuel-cell system takes the fuel costs and includes the efficiency for fuel to generate electric power. When fuel is transformed into power, the cost function considers the efficiency of fuel cell. Fuel-cell is the most efficient system among all fossil-fuel energy sources such as CHP units as:

$$f_{CHP}(P_{CHP}) = (C_{gasCHP} * \frac{P_{CHP}}{\mu_{CHP}}) - (C_{th}(\theta_{H/E})_{CHP} * P_{CHP}) \quad (10)$$

*b) O&M Cost Function:* The operation and maintenance (O&M) cost of the *n*th unit as a function of the power generated can be obtained by:

$$OM_n(p_n) = k_{OM} * p_n \quad (11)$$

*c) Pollution Function:* The emission function can be presented as the sum of all types of emissions considered, such as NOx, SOx, thermal emission, etc. In the present study, two important types of emission gases are taken into account. The amount of NOx and Sox emission is given as a function of generator output that is the sum of a quadratic and exponential function as follows:

$$E(P_G) = \sum_{i=1}^{Ne}(\alpha_i + \beta_i P_{G_i} + \gamma_i P_{G_i}^2) + \zeta_i exp(\lambda_i P_{G_i}) \quad (12)$$

*d) Power losses cost:* Hourly power losses considering network distribution losses cost are written as below:

$$PLC(t) = ((\sum_{b=1}^{Nb}[(P_b^2 + Q_b^2) \times Rb / V_b^2]) * C_{P_{loss}}) \quad (13)$$

### B. Reliability evaluation

The ENS cost is mainly considered as the cost of unsupplied demands over a time period (mainly one year). This cost is calculated as the amount of unsupplied energy multiplied by the forfeit as (16). Where, ENS specifies the amount of unsupplied energy (kWh) at stage t and AENS signifies the value of forfeit ($/kWh) for the unsupplied energy at T indicates set of hours at one day. It is clear that F2 indicates the daily cost in $/day. Then the average energy not supply or loss of load expectation (AENS or LOLE) is given by the following:

$$AENS = \sum_{1}^{N_c} ENS \times prob_i \quad (14)$$

The reliability of the network or energy index of reliability (EIR) is then given by the following:

$$EIR = 1 - \frac{AENS}{P_l} \quad (15)$$

$$C_{AENS} = \sum_{s=1}^{S}\sum_{b=1}^{Nb} Cint_b^s \lambda_b^s L_b (\sum_{res=1}^{N_{res}} P_{res}^s T_{res}^s + \sum_{rep=1}^{N_{rep}} P_{rep}^s T_{rep}^s) \quad (16)$$

It is worth remarking that many reliability indexes are defined in the distribution networks. But, in this paper, ENS is considered. Since, the objective function of this paper (Equation (16)) comprises several terms and all terms are given in "$/day". Therefore, the reliability cost should also be given in "$/day". As a result, ENS (kWh) multiplied by the value of forfeit ($/kWh) is equal to the reliability cost (ENS cost) in "$/h" or "$/day". In case of applying the other reliability indexes that their units are not $/day, the problem has to be reformulate as a multi-objective problem that allows different units. For instances, many electric distribution company use SAIFI or SAIDI because of end-users. In such case, the problem is reformulated as a multi-objective problem, since units of SAIFI and SAIDI are not $/year. We use the following equation to calculate the total interruption cost in MGs for the reliability index on a consumer side:

$$IC_{iday} = H_c \times C_{AENS} \quad (17)$$

$$F2 = \sum_{MG}(IC_{dey,MG}) \quad (18)$$

### C. Objective Function

As a result, the economic and reliability standpoints which have been explained in previous section are considered in two unique objective functions and converted into single objective model by weighted sum method as formulated below:

$$Z = [H_1 * F1 + H_2 * F2] \quad (19)$$

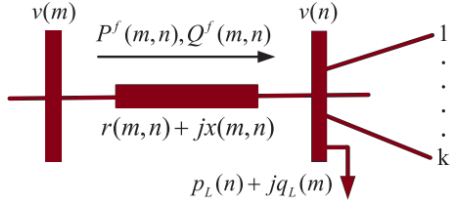

Fig. 2. Single line diagram of a radial distribution network.

*D. Constraints*

*a) Generation Capacity Constraint*: For stable operation, the real power output of each DG and ESS is restricted by lower and upper limits as follows:

$$P_{G_i}^{min} \leq P_{G_i} \leq P_{G_i}^{max}, \quad i=1,...N \quad (20)$$

$$SOC_{G_i}^{min} \leq SOC_{G_i} \leq SOC_{G_i}^{max}, \quad i=1,...M \quad (21)$$

*b) Power Balance Constraint*: The total electric power generation must cover the total electric power demand $P_D$ and the real power loss in distribution lines $P_{loss}$, hence:

$$\sum_{i=1}^{n} P_{G_i} - P_D - P_{loss} = 0 \quad (22)$$

Consider an electrical network as shown in Fig. 2, calculation of implies solving the load flow problem, which has equality constraints on real and reactive power at each bus as follows:

$$P_{i+1} = P_i - r_i \frac{(P_i^2 + Q_i^2)}{V_i^2} - P_{i+1} \quad (23)$$

$$Q_{i+1} = Q_i - x_i \frac{(P_i^2 + Q_i^2)}{V_i^2} - Q_{i+1} \quad (24)$$

In the above equations, we assume $P_i$ is generated by both RES-based DG units which are subject to uncertainties and controllable DG units, $Q_i$ is generated by controllable DG units.

## V. SOLUTION ALGORITHM

Cuckoo optimization algorithm is a meta-heuristic algorithm developed by Rajabioun in 2011 [20]. The basic idea of this algorithm is based on the obligate brood parasitic behavior of some cuckoo species in combination with the Levy flight behavior of some birds and fruit flies. Cuckoos are fascinating birds, not only because of the beautiful sounds they can make, but also because of their aggressive reproduction strategy. Some species such as the ani and guira cuckoos lay their eggs in communal nests, though they may remove others' eggs to increase the hatching probability of their own eggs. Quite a number of species engage the obligate brood parasitism by laying their eggs in the nests of other host birds. Some host birds can engage direct conflict with the intruding cuckoos. If a host bird discovers the eggs are not its own, it will either throw these alien eggs away or simply abandon its nest and build a new nest elsewhere. Each cuckoo starts laying eggs randomly in some other host birds' nests within her ELR. Fig. 3 gives a clear view of this concept. [21].

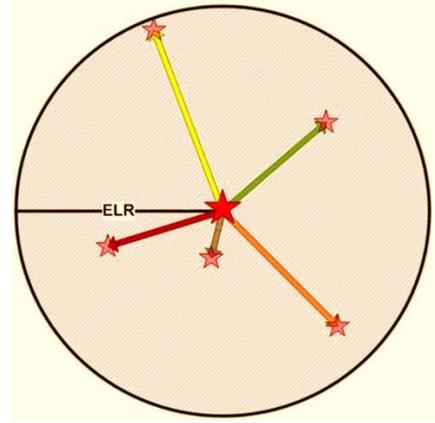

Fig. 3. Random egg laying in Egg Laying Radius (ELR).

## VI. RESULTS AND DISCUSSIONS

The PG&E 69-bus distribution network with three kinds of DGs units and ESSs is utilized for implementing the proposed stochastic scenario based model for optimal daily operation of MGs. Details about the PG&E 69-bus test system can be found in reference [22]. In this paper a typical low voltage (LV) MG portrayed in Figs. 4 is considered as the test system. The MG consists of different DER units such as the CHP, PV, WT and also NiMH-Battery based ESS. The system data is adopted from Ref. [23]. It is supposed that all DG units produce active power at unity power factor, neither requesting nor producing reactive power. Furthermore, there is a power exchange link between the mentioned MG and the utility (LV network) in order to trade energy during a day based on decisions of the microgrid central controller (MGCC). Table I offers the technical specifications of the utilized DER units including PV and WT in the mentioned MG. As shown Table I, although the PV and WT units don't use any fuel, their price is much higher than the other units. This fact is because of their high capital cost. The price of these units considers high to assign payback cost for the initial outlay or as maintenance and renewal costs. The hourly forecasted load demand inner the MG, the normalized forecasted output power of WT and PV for a typical day are shown in Fig. 5.

Furthermore, there is a power exchange link between the utility and the MG during the time step in the study period based on the decisions made by the MGCC. It is assumed that the MGCC purchases the maximum available power of the WT and PV at each hour of the day. Moreover, in order to make the analysis simpler, it is supposed that all the units work in electricity mode and no heat load demand is needed. It should be pointed out that all evolutionary methods require tuning of different algorithm parameters for their proper searching. A small change in these parameters may result in a large change in the algorithm performance. COA overcomes such difficulties, because it does not require any parameter for tuning. It means this algorithm reaches the optimal solution without adjusting any parameter.

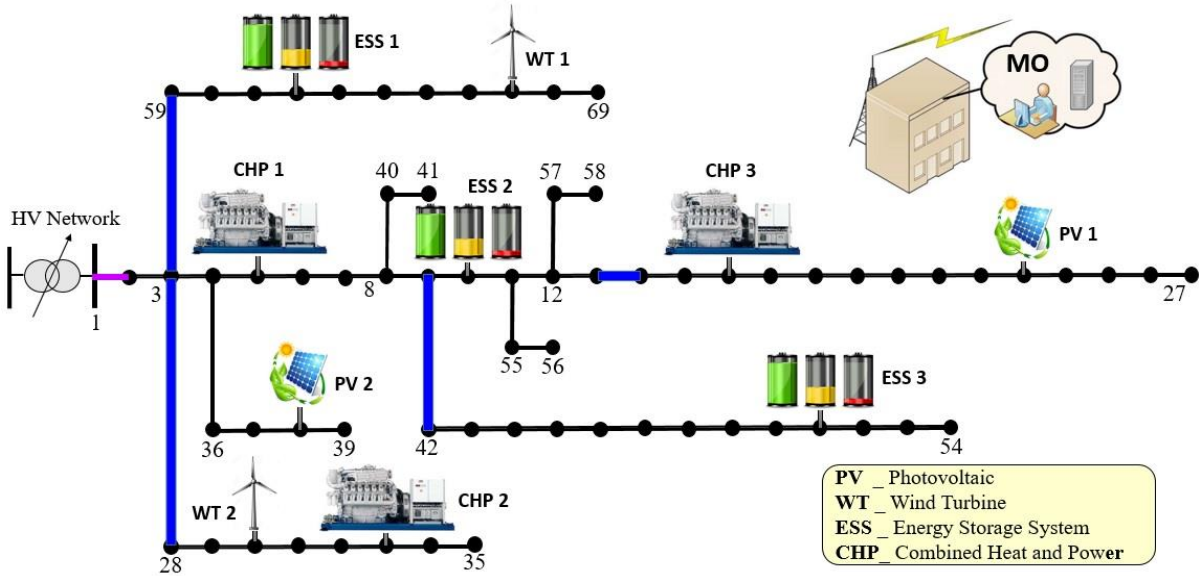

Fig. 4. Test smart MGs including DER units.

TABLE I. VALUES OF PARAMETERS IN WT AND PV.

| WT | Pr=250KW | $V_{CI}$=2m/s | $V_r$=14m/s | $V_{CO}$ = 25m/s |
|---|---|---|---|---|
|  | n=4 | α=3 | β=12 | - |
| PV | $P_{STC}$=250KW | $G_{STC}$=1000W/m² | K=0.001 | $T_c$=25 C |

In order to evaluate the performance of the extracted model, the peak load curve of the Iranian power grid on 28/08/2007 (annual peak load), Fig. 4, has been used for our simulation studies [7]. The load curve is divided into three different periods, namely low load period (00:00 am–8:00 am), off-peak period (8:00 am–6:00 pm) and peak period (6:00 pm–12:00 pm). The average value of generated power of utilized DER units consists of PV, WT, CHP and ESS units are shown in Fig. 5 for 24 hours.

From the point of MG operation, the stochastic approach by consideration of the intermittent nature of system components and loads can provide a more accurate solution for determining the allowance and optimum costs and emissions. The stochastic single objective MGs operation is solved by usage of the scenario generation scheme described in previous Section 2. 1000 24-h scenarios [22] are generated and subsequently trimmed to six batches. The number of scenarios is usually selected in a way that the coefficient of variation becomes small i.e. between 0.1% and 1% [23].

The voltage profile based on mean values has been shown in Fig. 6. According to this figure, with using the proposed method for daily energy management of MG, voltage profile of system is improved. Fig. 7 shows the power losses and ENS for proposed method and traditional energy management for smart distribution network. According to this figure, it is observed that by applying the proposed method for energy management the power losses of system are decreased. Based on correlation between network input and output variables, the output of the network parameters shows probabilistic behavior. It should be noted that the main disadvantage of the batches with enormous scenarios is the great CPU computational time requirement. With the help of fast and efficient computing tools i.e. scenario reduction technique and the COA algorithm, a huge stochastic model can be easily solved.

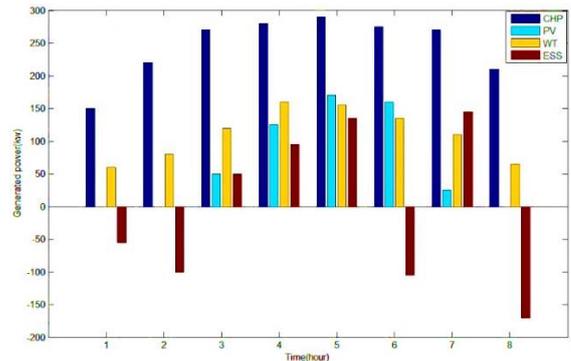

Fig. 5. Mean value of generated power by DER units in DEM.

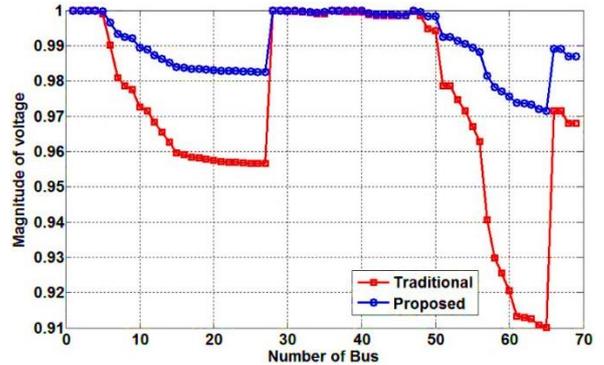

Fig. 6. Voltage profile of network based on mean value.

TABLE II. STATISTICAL ANALYSIS OF GENERATED POWERS.

| DER | PV | | WT | | CHP | | | ESS | | |
|---|---|---|---|---|---|---|---|---|---|---|
| Period (3h) | 1 | 2 | 1 | 2 | 1 | 2 | 3 | 1 | 2 | 3 |
| 1 | 0 | 0 | 26 | 34 | 52 | 43 | 55 | -15 | -23 | -18 |
| 2 | 0 | 0 | 33 | 47 | 73 | 68 | 80 | -27 | -39 | -35 |
| 3 | 31 | 19 | 60 | 60 | 78 | 95 | 93 | 22 | 15 | 13 |
| 4 | 68 | 55 | 85 | 78 | 85 | 96 | 98 | 27 | 38 | 35 |
| 5 | 89 | 76 | 66 | 89 | 88 | 102 | 100 | 34 | 42 | 58 |
| 6 | 84 | 78 | 53 | 59 | 94 | 95 | 87 | -29 | -33 | -37 |
| 7 | 12 | 11 | 44 | 67 | 79 | 104 | 92 | 54 | 45 | 42 |
| 8 | 0 | 0 | 23 | 38 | 66 | 75 | 65 | -58 | -63 | -49 |

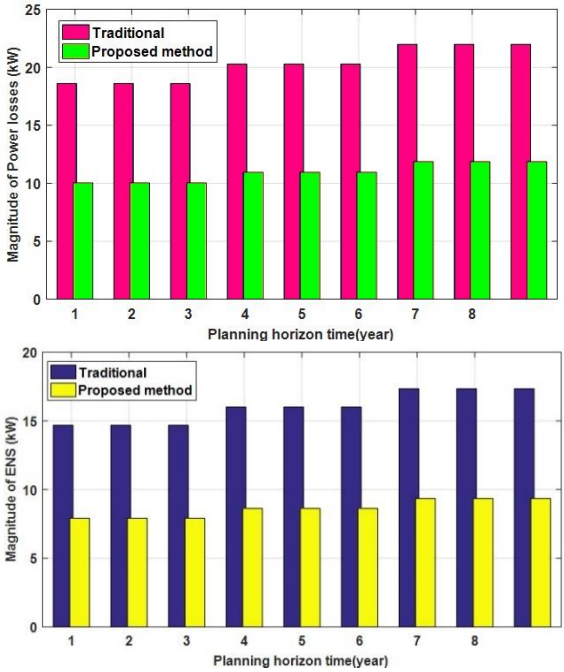

Fig. 7. Power losses and ENS of MG based on expected value.

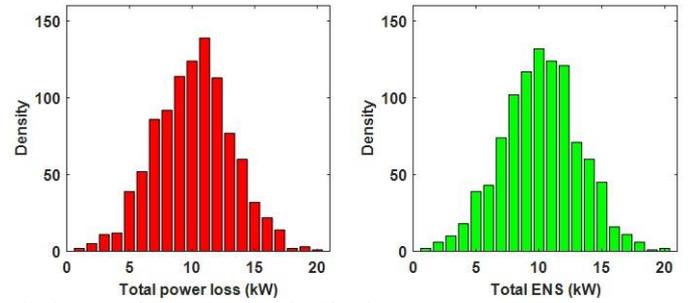

Fig. 8. PDFs of power losses and ENS value.

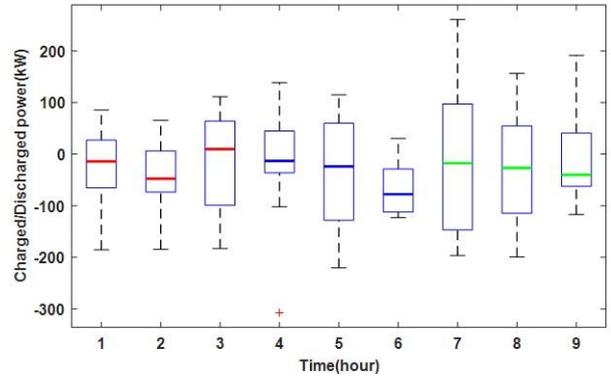

Fig. 9. Plot box of ESS over entire operation times.

In this paper, the stochastic approach tries to find expected optimal solutions for operation of MGs. The stochastic solutions may not be global optima solutions to the individual scenarios but they are a robust and also located near global solutions which this provides possible realizations of the uncertainties. In the optimization process, each scenario considers different values for load, WT and PV powers output and market price. Also, each scenario has different decision variables but the same state variables. The expected operation costs and emissions of these 30 scenarios are the output random variables which should be optimized simultaneously in the algorithm procedure. The power loss and ENS cost of MG can be shown in PDF form in Fig. 8. Fig. 9 illustrate the optimal charging/discharging of ESS units during operation horizon. It can be seen that in the off-peak times the ESS units charge and then in peak hours discharge into the network. Therefore, by optimal arbitrage of ESS units can be achived more profit for

MG operatore and it leads to more reduction in ENS value in peak periods. Also, it can be inferred from the results that modeling of the system uncertainties will increase the operation cost and emission values because stochastic procedure considers different most probable scenarios instead of one scenario (as the deterministic scheme). Indeed, approaching to the real conditions of the power system in the ante-scheduling studies will cost some expenses which are expectable. In other words, stochastic dynamic operation will concurrently consider the most probable scenarios. Besides, using the proposed stochastic framework, all 30 accepted scenarios according to their probability values contribute into the output random variable results, whereas the deterministic method relies on only one scenario. The 30 accepted scenarios capture more of the uncertainty spectrum of the power system, which is approximately four times more than that of the deterministic framework. So, the results of the stochastic framework are more realistic than the deterministic framework results.

## VII. Conclusion

In this paper, COA based method has been developed and applied to solve the dynamic operational scheduling of MGs with clean sources such as WT, PV, and ESS units for economic, reliable and emission operation problem by considering uncertainties including WT and PV powers output, and load demand over the 24 h study horizon. Also for getting closer to real condition as well as the reliability reasons, a reserve constraint is taken into account. The best advantage of the proposed algorithm is quick transfer of the information between agents which this gives more ability to the proposed algorithm in finding the global optima irrespective of the complexity of the problem. Besides, an innovative parameter setting technique is complemented to the original COA to cope with the drawback of premature convergence. This modification includes two powerful knowledge interaction strategies. Each individual according to a probability model chooses one of these methods to improve its knowledge. Since some RESs such as WT, PV have intermittent characteristic, approaches to analyze MGs would be stochastic rather than deterministic. The proposed approach shows how the proposed formulation works in comparison with an unreal-case-based deterministic technique. To take the uncertainties into account, a stochastic scenario-based method according to the Monte Carlo technique is implemented. Using the proposed model, possible scenarios of power system operating states are generated and a probability is assigned to each scenario. For a tradeoff between computation time and accuracy, a backward scenario reduction technique is utilized.

There are three possible avenues for future work arising from this paper, namely: 1) Considering other uncertainty modeling approaches such as IGDT, PEM and ROA; 2) Taking into account anther MGs management options like capacitor switching and network reconfiguration 3) integrating of PHEV parking lots in MGs as backup source to cope with operating uncertainties and 4) Using Multi objective optimization to create a tradeoff between various objectives.